\documentclass[leqno]{article}
\usepackage{amsmath,amsthm,amssymb,amscd}

\theoremstyle{plain}
\newtheorem{thm}{Theorem}[section]
\newtheorem{theorem}[thm]{Theorem}

\newtheorem{definition}[thm]{Definition}

\newtheorem{remark}[thm]{Remark}

\theoremstyle{defintion}

\def\dim{{\rm dim\ }}

\def\supp{\mathrm{supp}}
\def\n{\mathfrak{n}}
\def\h{\mathfrak{h}}
\def\z{\mathfrak{z}}

\def\NN{{\mathbb N}}

\def\RR{{\mathbb R}}
\def\CC{{\mathbb C}}
\def\HH{{\mathbb H}}

\def\Y{\mathcal{Y}}
\def\conv{\mathrm{conv}}

\def\B{\mathrm{B}}
\def\C#1{C_c\left(#1\right)}
\def\so{\mathfrak{so}}
\def\su{\mathfrak{su}}
\def\Ex{{\exp^\n}}
\def\msp{\medskip\par\noindent}
\def\hx{h_\text{max}}
\def\XX#1{\mathcal{X}#1}

\begin{document}

\title{The X-ray transform on $2$-step nilpotent Lie groups of higher rank}

\author{Norbert Peyerimhoff - Evangelia Samiou}




\maketitle

Dedicated to the memory of Sergio Console.

\begin{abstract}
We prove injectivity and a support theorem for the X-ray transform on $2$-step nilpotent Lie groups with many totally geodesic $2$-dimensional flats. The result follows from a general reduction principle for manifolds with uniformly escaping geodesics.
\end{abstract}


\section{Background}
\label{sec:geninform}

The X-ray transform of a sufficiently rapidly decreasing continuous
function $f$ on the Euclidean plane $\RR^2$ is a function $\XX{f}$
defined on the set of all straight lines via integration along these lines. More
precisely, if $\xi$ is a straight line, given by a point $x \in \xi$
and a unit vector $\theta \in \RR^2$ such that $\xi = x + \RR \theta$,
then
$$\XX{f}(\xi) =\XX{f}(x,\theta) = \int_{-\infty}^\infty f(x+s \theta)\, ds. $$
It is natural to ask about injectivity of this transform and, if yes,
for an explicit inversion formula. If $f(x) = O(|x|^{-(2+\epsilon)})$
for some $\epsilon > 0$, the function $f$ can be recovered via the
following inversion formula, going back to J. Radon \cite{Ra1917} in
1917:
\begin{equation} \label{eq:inv-form}
  f(x) = - \frac{1}{\pi} \int_0^\infty \frac{F_x'(t)}{t}\, dt,
\end{equation}
where $F_x(t)$ is the mean value of $\XX{f}(\xi)$ over all lines $\xi$
at distance $t$ from $x$:
$$ F_x(t) = \frac{1}{2\pi} \int_{S^1} \XX f(x+t \theta^\perp,\theta)\, d\theta, $$
where $(x,y)^\perp = (y,-x)$. Zalcman \cite{Za1982} gave an example of
a non-trivial function $f \in C^\infty(\RR^2)$ with
$f(x) = O(|x|^{-2})$ and $\XX{f}(\xi) = 0$ for all lines
$\xi \subset \RR^2$ and, therefore, the decay condition for the
inversion formula is optimal.

Under stronger decay conditions, it is possible to prove the following
support theorem (see \cite[Thm. 2.1]{Hel1965} or
\cite[Thm. I.2.6]{Hel2011}):

\begin{theorem}[Support Theorem] \label{thm:suppclass} Let $R > 0$ and
  $f \in C(\RR^2)$ with $f(x) = O(|x|^{-k})$ for all $k\in\NN$. Assume
  that $\XX{f}(\xi) = 0$ for all lines $\xi$ with $d(\xi,0) > R$. Then
  we have $f(x) = 0$ for all $|x| > R$.
\end{theorem}

Again, the stronger decay condition is needed here by a counterexample
of D.J. Newman given in Weiss \cite{We1967} (see also
\cite[Rmk. I.2.9]{Hel2011}). The Euclidean X-ray
transform plays a prominent role in medical imaging techniques
like the CT and PET (see, e.g., \cite{Ku2014}).

The X-ray transform can naturally be generalized to other complete,
simply connected Riemannian manifolds, by replacing straight lines by
complete geodesics. Radon mentioned in \cite{Ra1917} that there is an
analogous inversion formula in the (real) hyperbolic plane $\HH^2$,
where the denominator in the integral of \eqref{eq:inv-form} has to be
replaced by $\sinh(t)$ (see also
\cite[Thm. III.1.12(ii)]{Hel2011}). There is also an analogue of the
support theorem for the hyperbolic space (see
\cite[Thm. III.1.6]{Hel2011}), valid for functions $f$ satisfying
$f(x) = O(e^{-k d(x_0,x)})$ for all $k \in \NN$ and $x_0 \in \HH^n$.

In the case of a continuous function $f$ on a {\em closed} Riemannian
manifold $X$, the domain of $\XX f$ is the set of all {\em closed}
geodesics. Continuous functions $f$ can only be recovered from their
X-ray transform $\XX f$ if the union of all closed geodesics is dense in
$X$. But this condition is not sufficient as the following simple
example of the two-sphere ${\mathbb S}^2$ shows. Every {\em even}
continuous function $f$ on ${\mathbb S}^2$ (i.e., $f(-x)=f(x)$) can be
recovered by its integrals over all great circles. This fact and a
solution similar to \eqref{eq:inv-form} goes back to Minkowski 1911
and Funk 1913 (see \cite[Section II.4.A]{Hel2011} and the references
therein). But, on the other hand, it is easy to see that $\XX f$ vanishes
for all {\em odd} functions, so the restriction to even functions is
essential. For injectivity and support theorems of the X-ray transform
on compact symmetric spaces $X$ other than ${\mathbb S}^n$ see, e.g.,
\cite[Section IV.1]{Hel2011}. Injectivity properties of the extended
X-ray transform for symmetric $k$-tensors on closed manifolds (with
respect to the solenoidal part) play an important role in connection
with {\em spectral rigidity} (see \cite{GK1980}) and were proved for
closed manifolds with Anosov geodesic flows (see \cite[Thms 1.1 and
1.3]{DSh2003} for $k=0,1$) or strictly negative curvature (see
\cite{CrSh1998} for arbitrary $k\in\NN$).

Another class of manifolds for which the X-ray transform and its
extension to symmetric $k$-tensors has been studied are {\em simple
  manifolds}, i.e., manifolds $X$ with strictly convex boundary and
without conjugate points (see \cite{Sh1994}). An application is the {\em
  boundary rigidity problem}, i.e., whether it is possible to
reconstruct the metric of $X$ (modulo isometries fixing the boundary)
from the knowledge of the distance function between points on the
boundary $\partial X$.  Solenoidal injectivity is known for $k=0,1$
for all simple manifolds (see \cite{Muh1977} and \cite{AnRo1997}), and
for all $k \in \NN$ for surfaces \cite{PSU2013} and for negatively
curved manifolds \cite{PS1988}. There are also support type theorem
for the X-ray transform on simple manifolds (see
\cite{krishnan,krishplam} and \cite{UV2015} and the references
therein). A very recommendable survey with a list of open problems is
\cite{PSU2014}.

\section{A reduction principle for manifolds with uniformly escaping geodesics}
\label{sec:redprinc}

In this article, we will only consider complete Riemannian manifolds
$X$ whose geodesics escape in the sense of e.g. \cite{Wojtkowski1982},
\cite{Wojtkowski1978}, \cite{knigla}, in a uniform way. Simply
connected manifolds without conjugate points have this property, but
we like to stress that the main examples in this article will be {\em
  manifolds with conjugate points}. Geodesics will always be
parametrized by arc length.

\begin{definition}
  A Riemannian manifold $X$ has {\em uniformly escaping geodesics} if
  for each $r\in\RR^+_0$ there is $P(r)\in\RR^+_0$ such that for every
  geodesic $\gamma\colon\RR\to X$ and every $t>P(r)$, we have
  $d(\gamma(t),\gamma(0))>r$. We call $P$ an {\em escape function} of $X$.
\end{definition}

The smallest such function $P$,
$$ P(r) := \sup \{t \ge 0 \mid \exists \text{ geodesic } \gamma\colon\RR\to X,
d(\gamma(0),\gamma(t))\leq r \}$$
is thus required to be finite for all $r$. After time $P(r)$ every
geodesic has left a closed ball $\B_r(p)$ of radius $r\in \RR^+_0$ around its
center $p \in X$. The function $P$ increases and satisfies
$P(r)\geq r$. Note that $P$ may not be continuous.

Manifolds with this property must be simply connected and
non-compact. As mentioned earlier, simply connected Riemannian
manifolds without conjugate points have this property with escape
function $P(r)=r$.

The class of compactly supported continuous functions on such a
manifold is preserved under restriction to totally geodesic immersed
submanifolds. Thus if $f$ is a compactly supported continuous function
on $X$, say $\supp(f)\subset\B_r(p)$ for some $p\in X$ and $r>0$, and
$\phi\colon Y\to X$ a totally geodesic isometric immersion, then $f$
has compact support on $Y$ and
$\supp(f\circ\phi)\subset \B^Y_{P(r)}(p)$.  In particular, this holds
for geodesics (as $1$-dimensional immersions) and the integral of $f$
over any geodesic in $X$ is thus defined.

Before we formulate the reduction principle, let us first fix some
notation. The unit tangent bundle of $X$ is denoted by $SX$. For a
Riemannian manifold $X$ let $\C{X}$ be the space of all continuous
functions $f\colon X\to\CC$ with compact support. By $G(X)$ we denote
the set of (unparametrized oriented) geodesics, i.e.
$$G(X)=\left\{\gamma(\RR)\mid \gamma\colon\RR\to X\text{ geodesic }\right\} $$
The X-ray transform of $f\in\C{X}$ is the function
$\XX{f}\colon G(X)\to\CC$ with
$$\XX{f}(L)=\int_L f=\int_{-\infty}^{+\infty}f(\gamma(t))dt$$
if $L=\gamma(\RR)$ and $\gamma$ a unit speed geodesic.  \msp

\begin{definition} Let $r_0 \ge 0$ and $\sigma\colon[r_0,\infty)\to\RR^+_0$
  be a function.  We say that {\em the $\sigma$-support theorem holds on
    $X$} if for $p\in X$ and $f\in\C{X}$, $r\in [r_0,\infty)$ we have that
  $\XX{f}|_{G(X\setminus\B_{\sigma(r)}(p))}=0$ implies
  $f|_{X\setminus\B_r(p)}=0$. We say that {\em $X$ has a support
    theorem} if this holds for a function $\sigma$ with
  $\lim_{r\to\infty}\sigma(r)=\infty$.
\end{definition}

\begin{remark} If $X$ has a $\sigma$-support theorem, then $X$ has a
  support theorem for all smaller functions as well. Moreover, we can
  always modify $\sigma: [r_0,\infty) \to \RR^+_0$ to be monotone
  non-decreasing. If $r_0 = 0$, i.e., $\sigma: \RR^+_0 \to \RR^+_0$,
  the $\sigma$-support theorem implies injectivity of the X-ray
  transform.
\end{remark}

Then we have the following reduction principle.

\begin{theorem}\label{thm:red-princ} Let $X$ be a complete, Riemannian
  manifold which has uniformly escaping geodesics with escape function $P$.
  \begin{enumerate}
  \item[(i)] Assume there exists, for every $x \in X$, a closed totally
    geodesic immersed submanifold $Y \subset X$ through $x$ such that
    the X-ray transform on $Y$ is injective. Then the X-ray transform
    on $X$ is also injective.
  \item[(ii)] Let $\mu\colon[r_0,\infty)\to\RR^+_0$ be a function with
    $\mu \ge P(0)$. Assume there exists, for every $v \in SX$, a
    closed totally geodesic immersed submanifold $Y \subset X$ with
    $v \in SY$ such that the $\mu$-support theorem holds on $Y$. Then
    a $\sigma$-support theorem holds on $X$ for any function
    $\sigma\colon[r_0,\infty) \to \RR^+_0$ with $P(\sigma(r)) \le \mu(r)$
    for all $r \ge r_0$. In particular, we can choose $\sigma$ to be
    unbounded if $\mu$ is unbounded.
  \end{enumerate}
\end{theorem}

\proof (i) is obviously true by restriction since all geodesics in Y
are also geodesics in $X$.

For (ii), let $f\in\C{X}$ and
$r\ge r_0$. We fix a point $p\in X$ and let $\Y_p$ be a set of closed
totally geodesic immersed submanifolds $Y$ with $\mu$-support theorem
and so that each geodesic through $p$ lies in one of the $Y\in\Y_p$.
\msp We then have
$$f|_{X\setminus\B^X_r(p)}=0$$
if
$$\forall\, Y\in\Y_p : f|_{Y\setminus\B^Y_r(p)}=0,$$
since, by assumption, each geodesic in $X$ is contained in some
$Y$. Now, by the $\mu$-support theorem in $Y\in\Y_p$, we have
$$f|_{Y\setminus\B^Y_r(p)}=0$$
if
$$\XX{f}|_{G(Y\setminus\B^Y_{\mu(r)}(p))}=0.$$
Since $X$ has uniformly escaping geodesics property, this is
guaranteed if
$$\XX{f}|_{G(X\setminus\B^X_{s}(p))}=0$$
for any $s \ge 0$ with $P(s)\leq\mu(r)$. Thus $X$ has a
$\sigma$-support theorem for any function $\sigma\colon[r_0,\infty) \to \RR^+_0$
satisfying $P(\sigma(r)) \le \mu(r)$.
\endproof

\begin{remark} \label{rmk:contescfunc-optsigma}
If the escape function $P\colon\RR^+_0\to\RR^+_0$ is left-continuous, i.e. $\lim_{s \nearrow r}P(s)=P(r)$, we can choose $\sigma(r)=\sup\{s\ge 0\mid P(s)\leq\mu(r)\}.$
\end{remark}

\section{Applications of the reduction principle}
In this section we demonstrate that many interesting examples can be
derived by the reduction principle from $\RR^2$ and $\HH^2$. The X-ray
transform on the euclidean and on the hyperbolic plane is injective
and both have a $\mu$-support theorem with $\mu(r)=r$. This follows
directly from the euclidean or hyperbolic version of Radon's classical
inversion formula \eqref{eq:inv-form}, or Theorem \ref{thm:suppclass}.
\msp
If $X=X_1\times X_2$ is the product of two Riemannian manifolds of positive dimension with uniformly escaping geodesics, with escape functions $P_1$ and $P_2$ respectively, then $X$ has uniformly escaping geodesics with function $P$ satisfying
\begin{multline*}
$$\max\{P_1(r),P_2(r)\}\leq P(r)=\sup\left\{\sqrt{P_1(r_1)^2+P_2(r_2)^2}\,\mid\,
  r_1^2+r_2^2=r^2\right\}\\ \leq P_1(r)+P_2(r).
\end{multline*}
Each vector $v \in S(X_1 \times X_2)$ lies in a $2$-flat $F \subset
X_1 \times
X_2$, i.e. a totally geodesic immersed flat submanifold. By the
reduction principle, the $\sigma$-support theorem holds on $X_1\times
X_2$ for any function $\sigma$ with $P(\sigma(r)) \leq r$ for all $r
\in
[P(0),\infty)$. Note that this result does {\em not} require that
there are support theorems for the X-ray transforms on the factors
$X_1$ and $X_2$.  \msp

The reduction principle can also be applied to symmetric spaces of
noncompact type. These spaces have no conjugate points and each of
their geodesics is contained in a flat of dimension at least $2$ if
their rank is at least $2$. In non-compact rank-$1$ symmetric spaces
each geodesic is contained in a real hyperbolic plane. Therefore, the
reduction principle yields injectivity of the X-ray transform and a
support theorem with $\sigma(r)=r$ (\cite{Hel1981}, also
\cite[Cor. IV.2.1]{Hel2011}). \msp

Another interesting family are noncompact harmonic manifolds, which do
not have conjugate points. Prominent examples in this family are
Damek-Ricci spaces. In \cite{Rou2010}, Rouviere used the fact that
each geodesic of a Damek-Ricci space is contained in a totally
geodesic complex hyperbolic plane $\CC {\mathbb H}^2$ to obtain a
support theorem with $\sigma(r)=r$ for Damek Ricci spaces.
\msp

The main result in this article is about injectivity of the X-ray transform and a support theorem for a certain class of $2$-step nilpotent Lie groups with a left invariant metric and higher rank introduced in \cite{Sam2002}. By\cite{OS1974}
these spaces have conjugate points. Therefore, the methods of \cite{krishnan} do not immediately apply to these spaces. The spaces in \cite{Sam2002} differ also significantly from Heisenberg-type groups which do not even infinitesimally have higher rank.

\subsection{$2$-step nilpotent Lie groups have uniformly escaping \\geodesics.}\label{2sn}

The Lie algebra of a $2$-step nilpotent Lie algebra $\n$ splits
orthogonally as $\n=\h\oplus\z$, $\z=[\n,\n]$ the commutator and
$\h=\z^\perp$ its orthogonal complement. We can thus view
$\z\subset\so(\h)$ as a vectorspace of skew symmetric endomorphisms of
$\h$. We have
$$\langle [h,k]\mid z\rangle = \langle zh\mid k\rangle$$
for $h,k\in\h$, $z\in\z$. We show that $2$-step nilpotent Lie groups
have uniformly escaping geodesics, hence the X-ray transform for all
functions with compact support is defined.
\begin{theorem}\label{thm:nilpgeod}
Let $N$ be a simply connected $2$-step nilpotent Lie group with Lie algebra $\n=\z\oplus\h$, $\z\subset\so(h)$. Then $N$ has uniformly escaping geodesics with a continuous escape function $P$.
\end{theorem}
\proof We will prove that for each $r\in\RR_0^+$ there is $P(r)\in\RR^+$ such that every geodesic $\gamma$ with $\gamma(0)=e$ (the neutral element of $N$) we have that $d(\gamma(t),e)\leq r$ implies $t\leq P(r)$.

We denote by $\Ex\colon\n\to N$ the exponential map of the Lie group. Since $N$ is simply connected nilpotent this is a diffeomorphism. In particular,\\ $(\Ex)^{-1}(\B_r(e))\subset\B^\n_{\rho(r)}(0)$ for some increasing continuous function $\rho(r)\colon\RR_0^+\to\RR_0^+$ with $\rho(0)=0$.
We will show that there is $P(r)$ such that for every geodesic $\gamma$ in $N$ with $\gamma(0)=e$, the curve $(\Ex)^{-1}\circ\gamma$ has left $\B^\n_{\rho(r)}(0)$ after time $P(r)$.
\par
From \cite{kaplan}
 for a geodesic $\gamma(t)=\Ex(z(t)+h(t))$ with $z(t)\in\z$, $h(t)\in\h$, $\gamma'(0)=z_0+h_0$, we have
$$h''(t)=z_0h'(t),$$
$$z'(t)=z_0+\frac{1}{2}[h(t),h'(t)],$$
which we need to solve subject to the initial conditions
$$\gamma(0)=\Ex(z(0)+h(0))=e\text{ hence }z(0)=0=h(0),$$
$$\gamma'(0)=z_0+h_0=z'(0)+h'(0),$$
so that $\|z_0\|^2+\|h_0\|^2=1$. The solution to the first equation is
$$h(t)=\left((e^{tz_0}-1) z_0^{-1}\right)h_0.$$
Note that this is well defined even if $z_0$ is not invertible. Inserting this into the second equation gives
$$z'(t)=z_0+\frac{1}{2}\left[\left((e^{tz_0}-1) z_0^{-1}\right)h_0, e^{tz_0}h_0\right].$$
Taking the scalar product of this with $z_0$ gives
$$\langle z'(t)\mid z_0\rangle=\|z_0\|^2+\frac{1}{2}\langle z_0\mid\left[\left((e^{tz_0}-1) z_0^{-1}\right)h_0, e^{tz_0}h_0\right]\rangle$$
$$=\|z_0\|^2+\frac{1}{2}\langle z_0(e^{tz_0}-1) z_0^{-1}h_0\mid e^{tz_0}h_0\rangle$$
$$=\|z_0\|^2+\frac{1}{2}\|h_0\|^2-\frac{1}{2}\langle h_0\mid e^{tz_0}h_0\rangle,$$
since $e^{tz_0}$ is orthogonal. In order to compute $\langle z(t)\mid z_0\rangle$, we integrate,
$$\langle z(t)\mid z_0\rangle=t\|z_0\|^2+\frac{t}{2}\|h_0\|^2+\frac{1}{2}\langle h_0\mid (1-e^{tz_0})z_0^{-1}h_0\rangle.$$
It follows that
$$z(t)=tz_0+\frac{t\|h_0\|^2+\langle h_0\mid (1-e^{tz_0})z_0^{-1}h_0\rangle}{2\|z_0\|^2}z_0+w(t)$$
with $w(t)\in\z$ perpendicular to $z_0$. Hence, in the norm $\|\cdot\|$ of $\n$, we can estimate
\begin{multline*}
\|z(t)+h(t)\|^2\geq\|\left((e^{tz_0}-1)z_0^{-1}\right)h_0\|^2+ \\ +\frac{1}{4\|z_0\|^2}\left(2\|z_0\|^2t+t\|h_0\|^2+\langle h_0\mid (1-e^{tz_0})z_0^{-1}h_0\rangle\right)^2.
\end{multline*}
We split $\h=\oplus_{\lambda\in\RR} E(z_0,i\lambda)$ into the eigenspaces of $z_0$ and let $\hx\in E(z_0,i\lambda)=\CC$ be the largest component of $h_0$, $i\lambda$ the corresponding eigenvalue. Thus $|\hx|^2\geq\frac{1}{\dim\h}\|h_0\|^2$. Disregarding all other components, we estimate
\begin{multline*}
\|z(t)+h(t)\|^2\geq\\ \geq \left|\frac{e^{it\lambda}-1}{i\lambda}\right|^2|\hx|^2+\frac{1}{4\|z_0\|^2}\left(2\|z_0\|^2t+t\|\hx|^2+\text{Re}\left(\frac{1-e^{it\lambda}}{i\lambda}\right)|\hx|^2\right)^2
\end{multline*}
$$=\frac{2-2\cos(\lambda t)}{\lambda^2}|\hx|^2+\frac{1}{4\|z_0\|^2}\left(2\|z_0\|^2t+\left(t-\frac{\sin(t\lambda)}{\lambda}\right)|\hx|^2\right)^2$$
$$=\|z_0\|^2t^2+|\hx|^2\left(\frac{2-2\cos(\lambda t)}{\lambda^2}+t\left(t-\frac{\sin(\lambda t)}{\lambda}\right)+\frac{|\hx|^2}{4 \|z_0\|^2}\left(t-\frac{\sin(\lambda t)}{\lambda}\right)^2\right).$$
We now consider the cases:
\msp
$\|z_0\|^2\geq\frac{1}{2}$: Then $\|z(t)+h(t)\|^2\geq\frac{1}{2}t^2$.
\msp
If $\|z_0\|^2\leq\frac{1}{2}$, then $\|h_0\|^2=1-\|z_0\|^2\geq\frac{1}{2}$, hence $|\hx|^2\geq\frac{1}{2\dim\h}$. We can therefore estimate
$$\|z(t)+h(t)\|^2\geq\frac{1}{2\dim\h}\left(\frac{2-2\cos(\lambda t)}{\lambda^2}+t\left(t-\frac{\sin(\lambda t)}{\lambda}\right)+\frac{1}{4\dim\h}\left(t-\frac{\sin(\lambda t)}{\lambda}\right)^2\right).$$
If $\lambda=0$ the bracket evaluates to $t^2$, hence $\|z(t)+h(t)\|^2\geq \frac{1}{2\dim\h}t^2$.
\msp
If $0\leq t\leq\frac{\pi}{2\lambda}$ then $\cos(\lambda t)\leq 1-\frac{1}{2}(\lambda t)^2$. The other two summands are always nonnegative. Hence in this case,
$$\|z(t)+h(t)\|^2\geq \frac{t^2}{2\dim\h}.$$
\msp
If $t>\frac{\pi}{2\lambda}$ then $t-\frac{\sin(\lambda t)}{\lambda}\geq\left(\frac{\pi}{2}-1\right)t$. Observing that the rightmost and the leftmost summand are nonnegative, we get in this case that
$$\|z(t)+h(t)\|^2\geq \frac{(\pi-2)t^2}{4\dim\h}.$$
\msp
Thus we have shown that
$$\|z(t)+h(t)\|^2\geq t^2\min\left\{\frac{1}{2},\frac{1}{2\dim\h},\frac{\pi-2}{4\dim\h}\right\}=t^2\frac{\pi-2}{4\dim\h}.$$
Thus the curve $(\Ex)^{-1}(\gamma(t))=z(t)+h(t)$ has left $\B^\n_{\rho(r)}(0)$ after time $t=P(r):=\rho(r)\sqrt{\frac{4\dim\h}{\pi-2}}$.
\endproof

\def\t{\mathfrak{t}}
\subsection{X-ray transform on certain $2$-step nilpotent Lie groups}
Let $\h=\RR^{2q}=\CC^q$ and $\z=\t_{q-1}\subset\su(q)\subset\so(2q)$ be the Lie algebra of the maximal torus of $SU(q)$ and consider the $2$-step nilpotent Lie group $N_q$ with Lie algebra $\n_q=\z\oplus\h=\t_q\oplus\RR^{2q}$ endowed with a left invariant metric. In \cite{Sam2002}, it was shown that for every $q\in\NN$, $q\geq 3$, the Lie group $N_q$ has the property that each geodesic is
contained in a totally geodesic immersed $2$-dimensional flat
submanifold. The reduction principle and Theorem \ref{thm:nilpgeod} yield
\begin{theorem} The X-ray transform on $N_q$ is injective and has a support theorem.
\end{theorem}

\begin{remark} \label{rem:sigmaexample}
  Since the escape function $P$, defined at the end of the
  proof of Theorem \ref{thm:nilpgeod}, is continuous, $N_q$ admits a
  $\sigma$-support theorem with
  $\sigma(r) = \sup\{s \ge 0 \mid P(s) \le r\}$, due to Remark
    \ref{rmk:contescfunc-optsigma}.

    Moreover, the $\sigma$-support theorem can be extended to {\em
      general compact sets} (not only metric balls) in $N_q$. This
    extension is based on the following direct consequence of the
    classical support theorem (Theorem \ref{thm:suppclass}) for the
    $X$-ray transform on $\RR^2$: Let $K_0 \subset \RR^2$ be a compact
    set and $\conv(K_0) \subset \RR^2$ be its convex hull. Let
    $f \in C(\RR^2)$ with decay conditions as in Theorem
    \ref{thm:suppclass}. Then
    $\XX f\vert_{G(\RR^2 \backslash K_0)} = 0$ implies
    $f\vert_{\RR^2 \backslash \conv(K_0)} = 0$ (see
    \cite[Cor. I.2.8]{Hel2011}). Using this fact, we conclude for any
    compact set $K \subset N_q$, any point $p \in N_q$, and any
    $f \in C_c(N_q)$ with $\XX f\vert_{G(N_q \backslash K)} = 0$ that
  $$ f|_{N_q \setminus {\rm conv}_p(K)} = 0, $$
  where
  $$ \conv_p(K) = \{ x \in X \mid 
  \text{$\forall\, Y \in \Y_p$ with $x \in Y$}: x \in \conv_Y(K \cap Y) \} $$
  and  
  \begin{itemize}
  \item $\Y_p$ is a set of totally geodesic immersions of submanifolds
  isometric to $\RR^2$ such that each geodesic through $p$ lies in one
  of the $Y \in \Y_p$,
  \item ${\rm conv}_Y(Z)$ denotes the convex hull
  of a subset $Z$ of $Y \cong \RR^2$. 
  \end{itemize}
  The proof is a straightforward modification of the proof of Theorem
  \ref{thm:red-princ}. The $\sigma$-support theorem is then just the
  special case $K = B_{\sigma(r)}(p)$, since then
  $\conv_p(K) \subset B_r(p)$.
\end{remark}

\bigskip

{\bf Acknowledgements:} We are grateful to G. Paternain and G. Knieper
for helpful comments and relevant references.

\bigskip

\begin{flushleft}

{\bf AMS Subject Classification: 53C20, 53C30, 22E25}\\[2ex]

%
Norbert Peyerimhoff\\
Department of Mathematical Sciences, Durham University, Science Laboratories\\
South Road, Durham, DH1 3LE, UK\\
e-mail: \texttt{norbert.peyerimhoff@durham.ac.uk}\\[2ex]

Evangelia Samiou\\
Department of Mathematics and Statistics, University of Cyprus\\
P.O. Box 20537, 1678 Nicosia, Cyprus\\
e-mail: \texttt{samiou@ucy.ac.cy}\\[2ex]

%

\end{flushleft}

\normalsize

\end{document}